\documentclass[review,3p,times]{elsarticle}

\usepackage{amssymb}

\usepackage[switch]{lineno}

\usepackage{hyperref}
\usepackage{graphicx}
\usepackage{dblfloatfix}
\usepackage{float}
\usepackage{graphicx}
\usepackage{times}
\usepackage{epsfig}
\usepackage{bm}
\usepackage{breakurl}
\usepackage{graphicx, subfigure}
\usepackage{amsmath,amssymb, amsfonts, euscript, mathrsfs}
\usepackage{subeqnarray}
\usepackage{cases}
\usepackage{algorithm}
\usepackage{algorithmic}
\usepackage{fancybox}
\usepackage[table]{xcolor}
\usepackage{physics}

\definecolor{lightgray}{gray}{0.9}
\definecolor{lightblue}{rgb}{0.98,0.98,1.0}

\usepackage[normalem]{ulem}

\newcommand{\executeiffilenewer}[3]{%
\ifnum\pdfstrcmp{\pdffilemoddate{#1}}%
{\pdffilemoddate{#2}}>0%
{\immediate\write18{#3}}\fi%
}

\graphicspath{ {./images/} }
\newcommand{\T}{\mathsf{T}}

\let\oldequation\equation
\let\oldendequation\endequation

\renewenvironment{equation}
  {\linenomathNonumbers\oldequation}
  {\oldendequation\endlinenomath}

\journal{}
\begin{document}

\begin{frontmatter}



\title{Regularization in Space-Time Topology Optimization for \\
Multi-Axis Additive Manufacturing}


 \author[Manchester]{Weiming Wang}
 \ead{wwmdlut@gmail.com}
 
  \author[PME]{Kai Wu}
 \ead{K.Wu-4@tudelft.nl}

 \author[PME]{Fred van Keulen}
 \ead{A.vanKeulen@tudelft.nl}

 \author[SDE]{Jun Wu \corref{mycorrespondingauthor}}
 \cortext[mycorrespondingauthor]{Corresponding author}
 \ead{J.Wu-1@tudelft.nl}

 \address[Manchester]{Department of Mechanical, Aerospace and Civil Engineering, The University of Manchester, Manchester, United Kingdom}
 \address[PME]{Department of Precision and
 Microsystems Engineering, Delft University of Technology, Delft, The Netherlands}
 \address[SDE]{Department of Sustainable Design Engineering, Delft University of Technology, Delft, The Netherlands}



\begin{abstract}
In additive manufacturing, the fabrication sequence has a large influence on the quality of manufactured components. While planning of the fabrication sequence is typically performed after the component has been designed, recent developments have demonstrated the possibility and benefits of simultaneous optimization of both the structural layout and the corresponding fabrication sequence. This is particularly relevant in multi-axis additive manufacturing, where rotational motion offers enhanced flexibility compared to planar fabrication. The simultaneous optimization approach, called space-time topology optimization, introduces a pseudo-time field to encode the manufacturing process order, alongside a pseudo-density field representing the structural layout. To comply with manufacturing principles, the pseudo-time field needs to be monotonic, i.e., free of local minima. However, explicitly formulated constraints proposed in prior work are not always effective, particularly for complex structural layouts that commonly result from topology optimization.

In this paper, we introduce a novel method to regularize the pseudo-time field in space-time topology optimization. We conceptualize the monotonic additive manufacturing process as a virtual heat conduction process starting from the surface upon which a component is constructed layer by layer. The virtual temperature field, which shall not be confused with the actual temperature field during manufacturing, serves as an analogy for encoding the fabrication sequence. In this new formulation, we use local virtual heat conductivity coefficients as optimization variables to steer the temperature field and, consequently, the fabrication sequence. The virtual temperature field is inherently free of local minima due to the physics it resembles. We numerically validate the effectiveness of this regularization in space-time topology optimization under process-dependent loads, including gravity and thermomechanical loads.

\end{abstract}


\end{frontmatter}


\section{Introduction}
\label{sec:intro}

The past few years have seen tremendous progress in computational design for additive manufacturing~\cite{Graziosi2023advancing} and topology optimization approaches in particular~\cite{Brackett2011SFF,Zegard2016SMO}. In topology optimization, structural design is formulated as an optimization problem of distributing a limited amount of material in a design domain to maximize its performance, for example, load-bearing capacity. The highly optimized designs often have complex geometry, which makes them difficult to produce. While the flexibility of additive manufacturing makes it possible to produce geometrically complex components, it is not without challenges. In the past decade, a focus in topology optimization has been on the incorporation of additive manufacturability, to avoid optimized designs that cannot be 3D printed directly or need costly post-processing. The additive manufacturability is commonly characterized by geometric features such as the overhang angle (e.g., ~\cite{Langelaar2016AM,Gaynor2016SMO,Qian2017IJNME,vandeVen2018SMO}) and enclosed voids (e.g.~\cite{Liu2015FME,Donoso2022CMAME}). Let us refer to the review article by Liu \textit{et al}.~\cite{Liu2018SMO} for developments up to 2018. 

A recent development in topology optimization for additive manufacturing has been on the integration of the physics involved in the layer-by-layer additive process into design optimization~\cite{Bayat2023PMS}. This means simulating the physics, such as self-weight and thermomechanical loads, acting on intermediate structures during the fabrication process. Allaire \textit{et al}.~\cite{Allaire2017JCP} proposed one of the first approaches in this emerging direction, addressing overhang constraints, which had typically been treated as a geometric criterion~\cite{Langelaar2016AM,Gaynor2016SMO,Qian2017IJNME,vandeVen2018SMO}. By restricting the deformation of intermediate structures due to their self-weight, the final optimized structure exhibits fewer overhanging features, reducing the need for auxiliary supports during fabrication. The self-weight of intermediate structures was considered similarly by Amir and Mass~\cite{Amir2018}, Bruggi \textit{et al}.~\cite{Bruggi2018FEAD}, and Haveroth \textit{et al}.~\cite{Haveroth2022CMAME}. Starting from self-weight, more sophisticated and computationally demanding physics in the fabrication process has been incorporated in topology optimization, particularly for mitigating undesired effects due to thermomechanical loads~\cite{Allaire2018MMMAC,Pellens2020SMO,Misiun2021CMAME,Miki2021FEAD}. These process-informed topology optimization approaches, by incorporating a model of the layer-by-layer manufacturing process in design optimization, provide a more detailed level of quality control compared to focusing on geometric (and topological) features of optimized designs.

In the aforementioned topology optimization approaches, the structural layout is optimized to generate high-performance components. In addition to structural design, the fabrication sequence, according to which the component is produced line-by-line and layer-by-layer, significantly influences manufacturability and the quality of fabricated components. In wire and arc additive manufacturing (WAAM), for example, substantial distortion may occur during fabrication as well as after cutting off the build plate~\cite{Williams2016MST}. The distortion is largely influenced by the fabrication sequence. Typically, the fabrication sequence is planned after the structural layout has been designed. This sequential workflow results in a gap between the quality of digital designs and their physical counterparts. Earlier work tackled this problem by assuming a prescribed fabrication sequence, e.g., considering factors like material anisotropy due to a prescribed toolpath strategy in the optimization of the structural layout~\cite{Liu2017CAD,Dapogny2019CMAME}. More recent efforts have started optimizing the toolpath for powder bed fusion~\cite{Boissier2020SMO,boissier2021concurrent} and layers for multi-axis additive manufacturing~\cite{Wang2020SMO,Wang2023CMAME}. Boissier \textit{et al}.~\cite{Boissier2020SMO,boissier2021concurrent} represented the toolpath as a continuous curve, formulating it as a shape optimization problem. 

Our work is motivated by emerging robot-assisted multi-axis additive manufacturing (e.g., WAAM), where the rotation of robotic arms offers more flexibility than traditional 3D printing, allowing not only translational but also rotational motion. With rotation, the fabrication sequence is not limited to planar layers. This brings in benefits such as turning an overhanging feature into a self-supporting one, by appropriately rotating the build plate~\cite{Dai2018ToG}. In our previous work~\cite{Wang2020SMO,Wang2023CMAME}, the fabrication sequence was encoded by a pseudo-time field, where a larger pseudo-time indicates that the corresponding location is to be materialized later in the fabrication process. Isocontours of the pseudo-time field segment the structural layout into (curved) layers. This parameterization is akin to density-based topology optimization in the sense that each location in the design domain carries an optimization variable, opening up the solution space of fabrication sequences. This parameterization also enables the simultaneous optimization of the structural layout and the fabrication sequence, i.e., by optimizing the density and time field concurrently, referred to as space-time topology optimization~\cite{Wang2020SMO}. 

\begin{figure*}[ht!]
\centering
\includegraphics[width=0.98\linewidth]{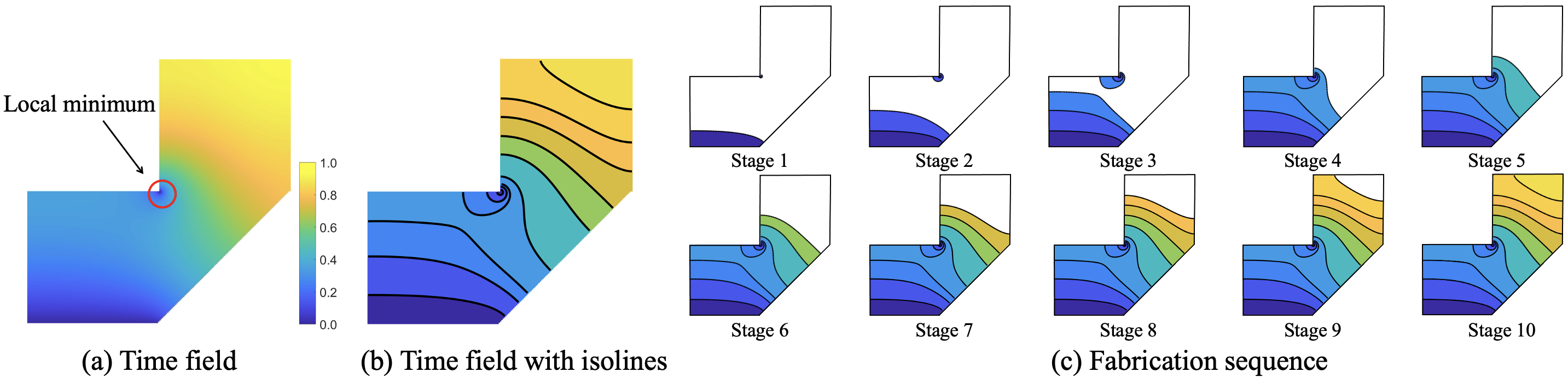}
\caption{(a) Illustration of a component and a pseudo-time field, with a smaller value indicating that the corresponding point is materialized earlier in fabrication. (b) The isolines of the pseudo-time field segment the part into a set of consecutive layers. (c) A series of intermediate structural layouts during the fabrication process. Note that, in stages 1 to 3, a patch of material, corresponding to a local minimum in the pseudo-time field, is in isolation from the rest of the component. This renders the fabrication sequence infeasible. }
\label{fig:problem}
\end{figure*}

Parameterizing the fabrication process by a pseudo-time field opens up a huge solution space. However, the sequence represented by the isocontours of the pseudo-time field is not guaranteed to be manufacturable. It is thus necessary to regularize the pseudo-time field to ensure that the optimized fabrication sequence meets manufacturing requirements. A critical requirement that we focus on in this paper is the continuity of fabrication: material can only be deposited on the already produced portion of the part, besides directly on the build plate. In the pseudo-time field, the violation of this requirement can be characterized by local minima, i.e., the surrounding points are associated with larger time values, and thus have not been materialized yet. A local minimum may present inside the domain or on the boundary of the domain. Both are undesirable, except it is adjacent to the build plate. An illustration of these violations is illustrated in Fig.~\ref{fig:problem}.

In our prior work~\cite{Wang2020SMO,Wang2023CMAME}, we addressed these violations by introducing explicitly formulated constraints. While these constraints function well in most cases, their effectiveness depends on the initialization of the pseudo-time field. In this paper, we present a novel regularization method to implicitly handle these constraints. Our starting point is to initiate the pseudo-time field using geodesic distances. Unlike Euclidean distances, geodesic distances are dependent on the structural layout and free of local minima. This initialization strategy works for fixed structural layouts, but is not directly applicable to space-time optimization, where the structural layout and the fabrication sequence are optimized concurrently. Because the structural layout evolves over iterations, and is not well defined during optimization, i.e., grey elements as encountered in density-based topology optimization. To solve these problems, we reformulate this initialization strategy using a partial differential equation, a heat equation in particular. This equation enables us to encode intermediate densities and evolving structural layouts in its parameters.

The remainder of this paper is organized as follows. In Section~\ref{sec:initialization}, we present a simplified problem setup of space-time optimization where the structural layout is fixed. This problem is solved by initialiating the fabrication sequence based on geodesic distances. In Section~\ref{sec:regularization}, this solution to the simplified problem setup is extended to address problems in space-time optimization where the structural layout and the fabrication sequence are optimized concurrently. The effectiveness of the proposed approach is evaluated in Section~\ref{sec:results} for design optimization in 2D and 3D under process-dependent self-weight and thermomechanical loads. In Section~\ref{sec:conclusions}, we conclude the paper with some ideas for future work.

\section{Fabrication sequence optimization}
\label{sec:initialization}

In this section, we consider a simplified version of space-time topology optimization. Here, the structural layout is prescribed, and we optimize the fabrication sequence to produce it. This reduced problem is intended to illustrate the problem of local minima, demonstrate the effects of the initial sequence in optimization, and highlight the efficacy of a novel initialization strategy. These insights inspire our approach to solving the problem in the full space-time topology optimization.

\subsection{Sequence parameterization} 
\label{subsec:definition}

We parameterize the fabrication sequence to produce a component using a pseudo-time field defined over the domain of the component. The component is discretized using a Cartesian grid, resulting in quadrilateral elements in 2D and hexahedral elements in 3D. To mimic the additive manufacturing process, each element $e$ is associated with a scalar value $t_e$, indicating the pseudo-time at which the material associated with the element is added to the structure. Thus, an element with a larger pseudo-time indicates that this element is materialized later than an element with a smaller one. For simplicity, we also refer to the pseudo-time field as the time field. 

\begin{figure}[tb]
\centering
\includegraphics[width=0.8\linewidth]{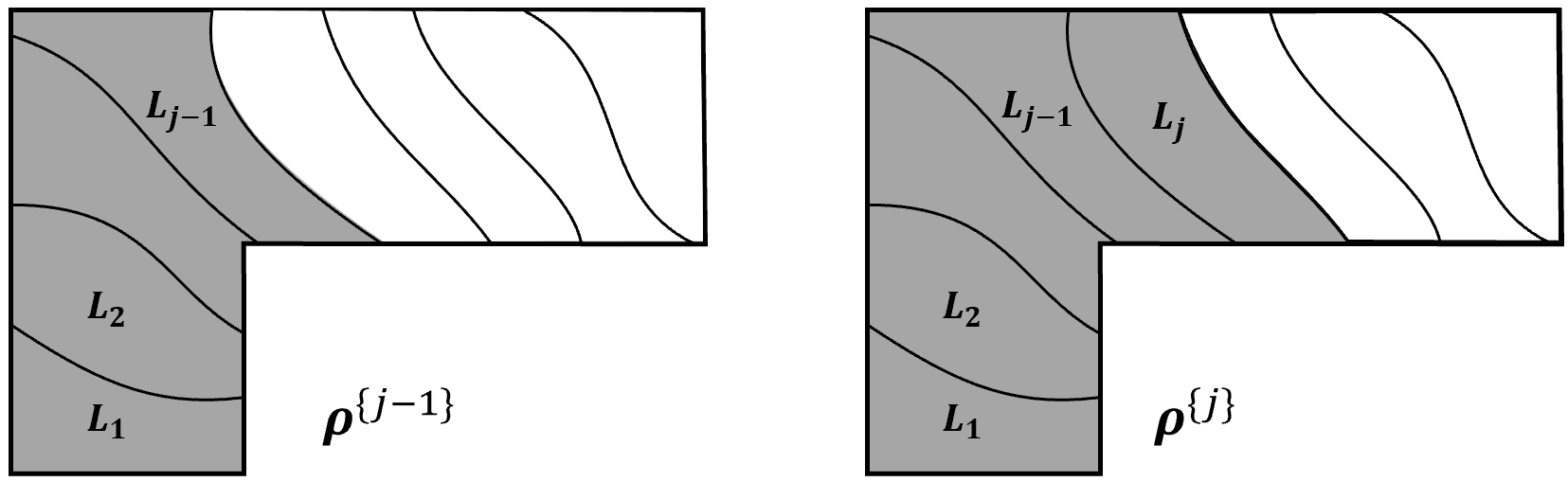}
\caption{Schematic of curved layers for fabricating an $L$-shaped component. Two successive intermediate structures during fabrication are illustrated on the left and right.}
\label{fig:layers}
\end{figure}

The time field $\bm{t}$ serves as design variables in sequence optimization. It is bounded to be between $0$ and $1$. The isocontours by $\eta_j = \frac{j}{N}$, with $j=0,...,N$, thus partition the component into $N$ layers. The partition is illustrated in Fig.~\ref{fig:layers}. Each layer consists of elements whose time value belongs to the corresponding time interval. For instance, the $j$-th layer is denoted as $L_j=\{e|t_e\in[\eta_{j-1},\eta_{j})\}, j=1,2,...,N$. We note that the time field is a continuous scalar field, mimicking the progressive additive manufacturing process. The isocontours of the time field, and the resulting layers between two consecutive isocontours, are not restricted to be planar but can be curved. 

We describe the intermediate structures during fabrication by a series of pseudo-density fields (or density fields, for simplicity), $\bm{\rho}^{\{j\}}$, with $j=0, 1, ..., N$. Two successive intermediate structures are illustrated in Fig.~\ref{fig:layers}. The superscript $\{j\}$ indicates the $j$-th stage during fabrication, i.e., layers from the first up to the $j$-th have been fabricated. $j=0$ refers to the start of fabrication, $\rho^{\{0\}}_e = 0$, $\forall e$, while $j=N$ refers to the completion of fabrication, $\rho^{\{N\}}_e = 1$,~$\forall e$. 

The evolving intermediate structure is derived from the time field. Let $\bm{1}$ denote the density field of a prescribed, solid component. Assuming that the density field and time field are discretized using the same Cartesian grid, the intermediate structure at the $j$-th stage is described by an entry-wise multiplication between the density and the time field,
\begin{equation}
    \bm{\rho}^{\{j\}} (\bm{t}) = \bm{1} \circ \bar{\bm{t}}^{\{\eta_j\}}, \quad j=1,...,N.
    \label{eq:rou_t}
\end{equation}
Here $\bar{\bm{t}}^{\{\eta\}}$ is a function of the time field $\bm{t}$. It converts time values less than $\eta$ to $1$ and time values greater than $\eta$ to $0$, in a differentiable manner. It is implemented using a smoothed Heaviside function, a commonly employed approach in density-based topology optimization~\cite{wang2011projection,wu2018infill}, 
\begin{equation}
\label{eq:projection}
    \bar{{t}}^{\{\eta\}} = 1 - \frac{\tanh(\beta_t \eta)+\tanh(\beta_t({t} - \eta))}{\tanh(\beta_t \eta)+\tanh(\beta_t({1} - \eta))},
\end{equation}
where $\beta_t$ is a positive number to control the projection sharpness. A continuation scheme is applied to increase this sharpness parameter after a certain number of iterations.

\subsection{Sequence optimization}
\label{subsec:sqeuence_formulation}

The sequence optimization is generically formulated as 
\begin{eqnarray}
\label{eq:sequence_optimization_framework}
 & \underset{\bm{t}}{\text{minimize}:} & c(\bm{t})=f(\bm{\rho}^{\{1\}}(\bm{t}), ..., \bm{\rho}^{\{j\}}(\bm{t}), ..., \bm{\rho}^{\{N\}}(\bm{t})), \\
 & \text{subject to}:     & t_e \geq \min \{t_i|i\in \mathcal{N}_e\}, \forall e \in \mathcal{M} \backslash \mathcal{S}_0, \label{eq:localminimum}\\
                        && t_e \leq \max \{t_i|i\in \mathcal{N}_e\}, \forall e \in \mathcal{M} \backslash \mathcal{B}, \label{eq:localmaximum}\\
                        && t_e \in [0,1], \forall e.
\end{eqnarray}
The objective concerns the evolving structure during fabrication $\bm{\rho}^{\{j\}}$, which is a function of the time field (cf. Eq.~\ref{eq:rou_t}). In Eq.~\ref{eq:localminimum} and Eq.~\ref{eq:localmaximum}, $\mathcal{N}_e$ denotes the set of elements adjacent to element $e$. Eq.~\ref{eq:localminimum} requests the time value of each element to be larger than the minimum among its adjacent elements. This requirement applies to all elements in the component ($\mathcal{M}$), except those designated as the start region for fabrication ($\mathcal{S}_0=\{e|t_e=0\}$), i.e., on the build plate. An element where this condition is violated indicates that the element needs to be materialized in midair, i.e., before any of its adjacent elements have been produced. 
Similarly, Eq.~\ref{eq:localmaximum} requests the time value of each element to be smaller than the maximum among its adjacent elements. This requirement applies to all elements in the component ($\mathcal{M}$), except those on the boundary of the component ($\mathcal{B}$). An element where this condition is violated indicates that the element needs to be materialized within an enclosed region, i.e., all of its adjacent elements have been produced.
Both Eq.~\ref{eq:localminimum} and Eq.~\ref{eq:localmaximum} are essential for the validity of a fabrication sequence. 

The challenges to solving this formulation are twofold. Firstly, the constraints (Eq.~\ref{eq:localminimum} and Eq.~\ref{eq:localmaximum}) involve $\min$ and $\max$ functions, which are non-differential. Secondly, the number of constraints is large, since both Eq.~\ref{eq:localminimum} and Eq.~\ref{eq:localmaximum} apply to all elements in the component. In previous work~\cite{Wang2020SMO,Wang2023CMAME}, these constraints were reformulated using approximation and aggregation schemes. For instance,
\begin{equation}
    \frac{1}{{n}(\mathcal{M})}\sum_{e\in \mathcal{M}}||t_e - \textrm{mean}_{i\in \mathcal{N}_e} (t_i)||^2 \leq \epsilon,
    \label{eq:meanConstraint}
\end{equation}
where $n$ denotes the number of elements in a set. The $\text{mean}$ function calculates the average of time values of adjacent elements, i.e., $\textrm{mean}_{i\in \mathcal{N}_e} (t_i)= \frac{1}{{n}(\mathcal{N}_e)} \sum_{i \in \mathcal{N}_e} t_i$.  
$\epsilon$ is a small constant for numerical stability (e.g., $\epsilon=10^{-6}$). While this approximation works well in most cases, its efficacy depends on the initialization of the fabrication sequence. 

\begin{figure*}[htb!]
\centering
\includegraphics[width=.96\linewidth]{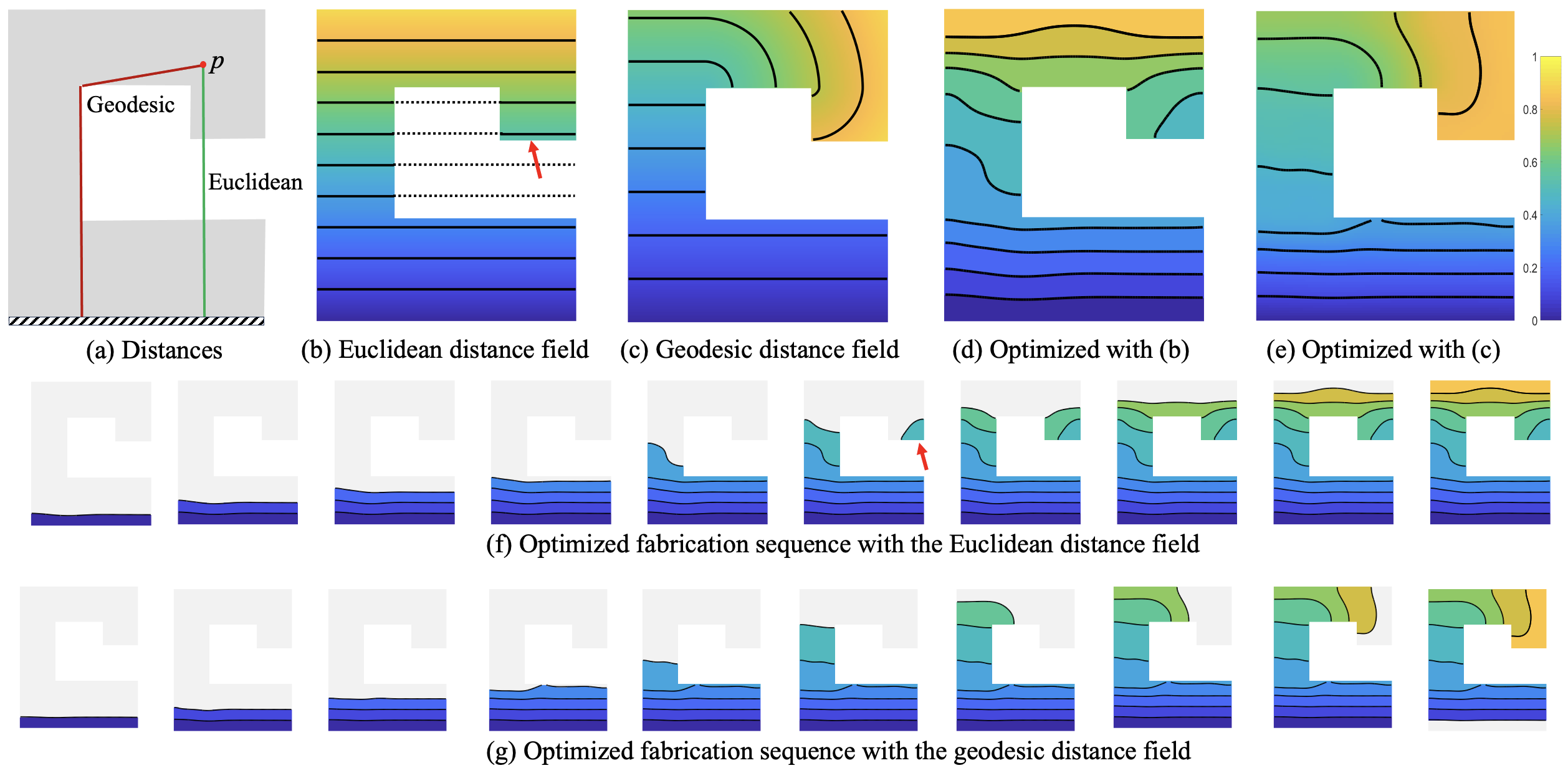}
\caption{Comparison of the Euclidean and the geodesic distance field as the initialization of the time field for fabrication sequence optimization. (a) A 2D component. The Euclidean distance and the geodesic distance from location $p$ to the build plate at the bottom are shown by the green line and the black polyline, respectively. (b) The Euclidean distance field. (c) The geodesic distance field. (d) The optimized time field using (b) as initialization. (e) The optimized time field using (c) as initialization. The sequences corresponding to the time fields in (d) and (e) are shown in (f) and (g), respectively.}
\label{fig:distance}
\end{figure*}

\subsection{Initialization of fabrication sequence}

The initialization of the fabrication sequence greatly influences the validity of optimized sequences, especially when the structural layout is complex. Distance fields serve as a good option for initialization as they can mimic the progression of additive manufacturing. A distance field is a scalar field, representing the minimal distance from each point in the domain to a certain shape, i.e., the build plate. Depending on the distance metric, the distance field has quite different features. Fig.~\ref{fig:distance} compares the effects of two distance metrics as initialization for sequence optimization. The component is illustrated in Fig.~\ref{fig:distance} (a), where a horizontal build plate is located at the bottom\footnote{In multi-axis additive manufacturing, the build plate is not necessarily flat nor located at the bottom.}. The two distance metrics under comparison are the Euclidean distance and geodesic distance. The Euclidean distance measures the shortest path (i.e., a straight line) between a point and the build plate. The Euclidean distance is irrespective of the shape of the component. In contrast, the geodesic distance measures the shortest path, within the component, from a point to the build plate. The shortest Euclidean path and geodesic path for a single point are illustrated in Fig.~\ref{fig:distance} (a). The corresponding distance fields, together with isocontours, are visualized in Fig.~\ref{fig:distance} (b) and (c). Both distance fields are normalized by dividing the distance values over the corresponding largest distance. In Fig.~\ref{fig:distance} (b), the isocontours of the Euclidean distance field are planar. On the right-hand side of the domain, there is an overhanging feature, and there the isocontours are separated by the void region. In contrast, in Fig.~\ref{fig:distance}~(c), the geodesic distance field is free of local minima. The isocontours of the geodesic distance field in the top-right of the component are curved. The optimized fabrication sequences starting from these two initial sequences are shown in Fig.~\ref{fig:distance} (d) and (e), while the corresponding layered fabrication processes are shown in (f) and (g), respectively. It can be seen in Fig.~\ref{fig:distance} (f) that during fabrication, an isolated layer appears at stage 6. This shows that the optimization is incapable of eliminating the local minimum inherited from the initial field. While this violates the constraint in Eq.~\ref{eq:localminimum}, it numerically satisfies the approximated version (Eq.~\ref{eq:meanConstraint}). Setting a smaller threshold ($\epsilon$) proves to be ineffective. Even though the approximated constraint is not precisely satisfied, the optimization is incapable of steering the location of the local minimum such that it disappears. This experiment shows the importance of a valid initial sequence that can mimic the progression of additive manufacturing and is free of local minima. The geodesic distance field fulfills these requirements.

\section{Regularization of space-time topology optimization with heat equations}
\label{sec:regularization}

As discussed in the previous section, using the geodesic distance field to initialize the time field is promising for optimizing the fabrication sequence when the structural layout is fixed. However, it is not directly applicable to space-time topology optimization, where the structural layout, including its topology, evolves during the optimization process. 

In this section, we present a novel approach to regularize the time field in order to avoid the occurrence of local minima in the context of space-time topology optimization. Specifically, we relate the time field to the solution of a heat equation, i.e., a temperature field. The temperature field is monotonic and free of local minima. Furthermore, the heat equation is defined on the structural layout, and the resulting temperature field is sensible to changes in the structural layout. The monotonicity and the sensibility to shape changes are essential features of the geodesic distance field. Furthermore, by adapting the parameters in the heat equation we can accommodate intermediate density values as encountered during optimization. Our adoption of the heat equation is inspired by advancements in computer graphics, where research has demonstrated that the heat equation offers a computationally efficient means for calculating an approximate geodesic distance field~\cite{2013Geodesics}. In topology optimization for additive manufacturing, heat equations were previously used for identifying enclosed voids and consequently avoiding them~\cite{VTM2015,VTM2020}. 

\subsection{Optimization workflow}
\label{subsec:workflow}

\begin{figure*}[htb!]
\centering
\includegraphics[width=0.98\linewidth]{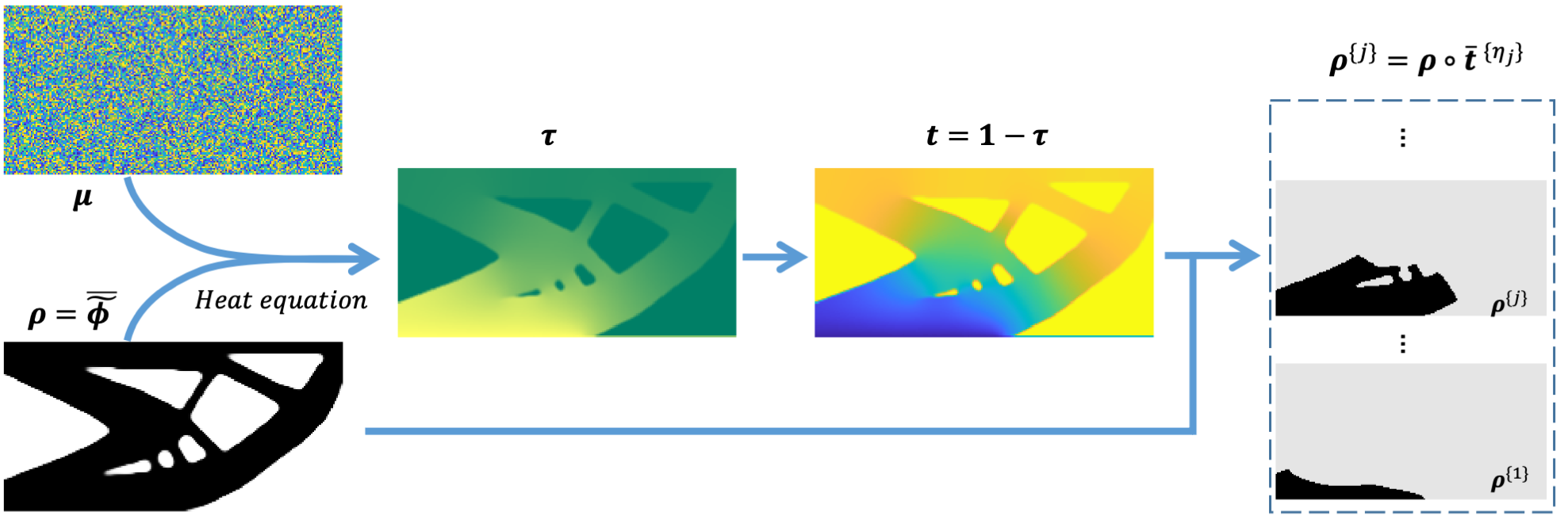}
\caption{The workflow of deriving a series of intermediate structures $\bm{\rho}^{\{j\}}$ from the thermal diffusivity field $\bm{\mu}$ and density field $\bm{\rho}$.}
\label{fig:workflow}
\end{figure*}

The workflow of space-time topology optimization using heat equations for regularizing the time field is illustrated in Fig.~\ref{fig:workflow}. We employ a density field ($\bm{\rho}$, $\rho_e\in[0,1]$) to represent the complete structure and a time field ($\bm{t}$, $t_e\in[0,1]$) to encode the fabrication sequence. While the density field eventually converges to either 0 or 1, the time field is continuous. The product of $\bm{\rho}$ and $\bm{t}$ is projected into a series of intermediate structures, using increasing threshold values $\eta_j = \frac{j}{N}$, with $j=0,...,N$ and $N$ being the prescribed number of layers. Rather than taking $\bm{\rho}$ and $\bm{t}$ as design variables, we introduce two auxiliary fields, $\bm{\phi}$ and $\bm{\mu}$. The density field $\bm{\rho}$ is computed from design field $\bm{\phi}$, following a density filter ($\tilde{\bm{\phi}}$) and a smoothed Heaviside projection ($\bar{\tilde{\bm{\phi}}}$). This process ($\bm{\phi} \to \tilde{\bm{\phi}} \to \bm{\rho}=\bar{\tilde{\bm{\phi}}}$) is commonly used in density-based topology optimization, and thus not included in the workflow illustration. 

The novel element in the workflow is the transformation from $\bm{\mu}$ and $\bm{\rho}$, through ${\bm{\tau}}$, to $\bm{t}$. We use a heat equation to obtain a monotonic time field defined on an evolving structural layout $\bm{\rho}$ with heterogeneous thermal diffusivity $\bm{\mu}$. The solution of the heat equation is a virtual temperature field $\bm{\tau}$. It is important to note that this virtual temperature field shall not be confused with the temperature field coming from the thermal process involved in metal additive manufacturing. Here the heat equation is meant to create a monotonic field, rather than predicting the actual temperature distribution during fabrication. Conceptually, it works as follows. An artificial heat source is placed on the build plate (at the bottom of the domain), from where the heat diffuses through the structure. This creates a continuous temperature distribution within the structure. The heat source has a constant temperature of $1$. At a location in the structure further away from the heat source, its temperature approaches $0$. The time field, ascending from the build plate, is defined as $\bm{t} = \bm{1} - \bm{\tau}$.

\subsection{Heat equation}
\label{subsec:heatequation}

To obtain a monotonic temperature field, we construct the following heat equation,
\begin{equation}
    \nabla \cdot (\kappa \nabla \tau) - \alpha_\texttt{T} \tau = 0.
\label{eq:heat_eq}
\end{equation}
Here $\nabla$ is the vector differential operator. This partial differential equation represents steady-state heat conduction within a heterogeneous medium, featuring spatially-varying thermal diffusivity $\kappa(x,y)$, and a drain term $\alpha_\texttt{T} \tau$, where $\alpha_\texttt{T}$ denotes a constant drain rate. The introduction of the drain term serves a crucial purpose. In the absence of this term, achieving a steady-state temperature distribution with spatial variation would necessitate employing an inhomogeneous Dirichlet boundary condition. However, when dealing with a structure that evolves during space-time topology optimization, there is no clear rationale for determining the ideal location of a heat sink. An alternative approach to achieve a spatially-varying temperature distribution is to resort to transient heat conduction, and calculate the temperature distribution at a specific time~\cite{2013Geodesics,Zhai2024CMAME}.

The heat equation is subject to the following boundary conditions: At the boundary adjacent to the build plate $\Gamma_{0}$, the temperature is constant, ${\tau\big|}_{\Gamma_{0}}=1$. Other boundary segments of the structure are thermally insulated, ${(\nabla \tau\cdot\overrightarrow{n})}\big|_{\Gamma\setminus\Gamma_{0}}=0$. Here $\nabla \tau$ represents the spatial gradient of the temperature field, and $\overrightarrow{n}$ denotes the normal direction on the boundary of the structure.

\begin{figure*}[htb!]
\centering
\includegraphics[width=0.9\linewidth]{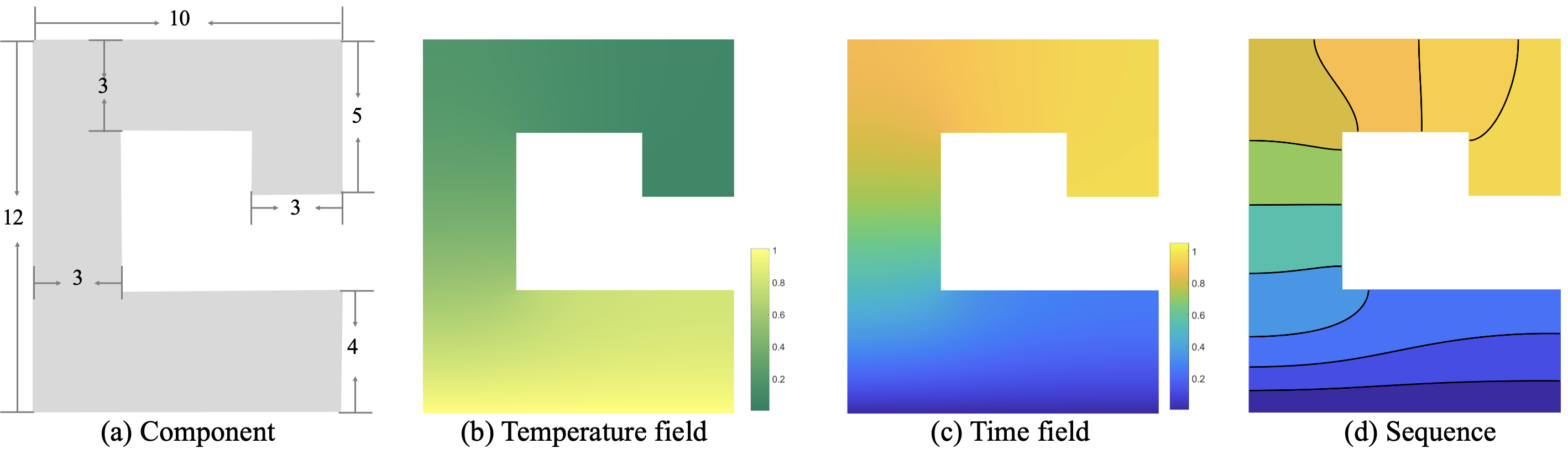}
\caption{Heat conduction on a given component (a). The start points are the nodes on the bottom boundary. The generated virtual temperature field by solving Equation~(\ref{eq:heatFEM}) is shown in (c). The light yellow indicates the largest temperature, and the dark blue indicates the smallest virtual temperature. The corresponding time field is shown in (c). The sequence generated with the time field shown in (d) which is composed of 10 layers.}
\label{fig:heat_sim}
\end{figure*}

The solution of the heat equation, and the corresponding time field together with fabrication sequences, are illustrated in Fig.~\ref{fig:heat_sim}. In this example the thermal diffusivity field ${\kappa}{(x,y)}$ is homogeneous. The temperature field and, consequently, the fabrication sequence are determined by the thermal diffusivity field. To allow the optimization of the fabrication sequence, the thermal diffusivity field is related to design variables $\bm{\mu}$, augmented by the unitless pseudo-density $\bm{\rho}$,
\begin{equation}
    \bm{\kappa}=\bm{\rho} \bm{\mu}.
\label{eq:heatconductivity}
\end{equation}
The multiplication by the pseudo-density field constrains heat conduction within the structural layout, excluding void regions from the fabrication sequence planning. Additionally, this allows pseudo-densities to take intermediate values during the optimization process, thereby facilitating gradient-based numerical optimization.

We note interesting similarities and large differences to the filters based on Helmholtz‐type partial differential equations (PDEs)~\cite{Lazarov2011PDEfilter}. In our formulation, the optimization variables (thermal diffusivity) in the PDE and the solution of the PDE (virtual temperature) have different physical interpretations. The desirable characteristic of the PDE is to attain a monotonic scalar field. In contrast, in the PDE filter~\cite{Lazarov2011PDEfilter}, the optimization variables and the solution both represent density fields. The intention behind it is to get a smoothed scalar field, avoiding checkerboard patterns.

\subsubsection{Numerical discretization}
\label{subsubsec:discretization}

We discretize Eq.~\ref{eq:heat_eq} using the finite element method, resulting in a linear form,
\begin{equation}
    \bm{K}_\texttt{T} \bm{T} = \bm{b},
\label{eq:heatFEM}
\end{equation}
where $\bm{T}$ is the nodal temperature vector. $\bm{b}$ is the thermal load vector corresponding to the boundary condition. $\bm{K}_\texttt{T}$ is the global system matrix, assembled from element matrices:
\begin{equation}
    \bm{K}_\texttt{T} = {\sum}_{e}{\bm{K}^{e}_\texttt{T}} +{\sum}_{e}{\bm{K}^{e}_\texttt{C}}.
\label{eq:heat_stiffness_matrix}
\end{equation}
Here the symbol ${\sum}$ stands for matrix assembly. The first element matrix corresponds to the thermal diffusivity. It is derived from the shape functions $N_i$ and $N_j$ by using the variational principle:
\begin{equation}
    \bm{K}^{e}_\texttt{T} = \int_{\Omega_e}{\kappa_e\left(\frac{\partial N_i}{\partial x}\frac{\partial N_j}{\partial x}+\frac{\partial N_i}{\partial y}\frac{\partial N_j}{\partial y}\right) d\Omega}.
\label{eq:heat_element_stiffness_matrix}
\end{equation}
The second element matrix, corresponding to the drain term, can be derived similarly,
\begin{equation}
    \bm{K}^{e}_\texttt{C} = \int_{\Omega_e}{\alpha_\texttt{T} N_i N_j d\Omega}.
\label{eq:heat_element_capacity_matrix}
\end{equation}
In our implementation, instead of calculating a consistent matrix $\bm{K}^{e}_\texttt{C}$ from the shape functions using Eq.~\ref{eq:heat_element_capacity_matrix}, we adopt a lumped matrix for simplicity, assuming that the drainage is lumped at the nodes. 

The nodal temperature distribution from solving Eq.~\ref{eq:heatFEM} is transformed into an elemental temperature field using bi-linear interpolation. In the matrix form, this is written as $\bm{\tau} = \bm{G} \bm{T}$. $\bm{G}$ is a sparse transformation matrix for converting nodal values to elemental values. Afterwards, the time field is computed by $\bm{t}=\bm{1}-\bm{\tau}$.

\subsubsection{Parameter analysis}
\label{subsubsec:parameterchoice}

\begin{figure}[tb]
\centering
\includegraphics[width=0.98\linewidth]{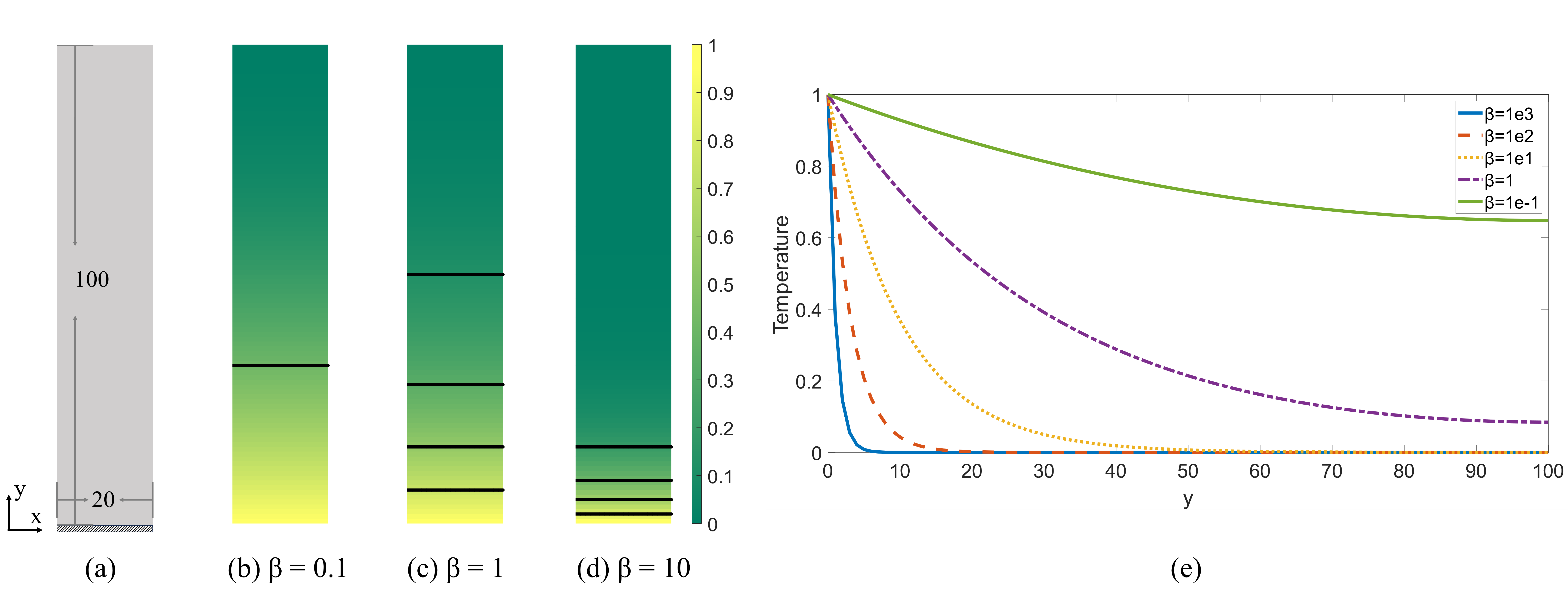}
\caption{The influence of the drain rate on the temperature distribution. (a) A rectangular component with a homogeneous thermal diffusivity. The virtual heat source with a fixed temperature of 1 is located at the bottom of the domain. (b, c, d) Virtual temperature distributions with $\beta = 0.1$, $1$, and $10$, respectively. The black lines represent the isocontours of the temperature fields at temperature values of 0.2, 0.4, 0.6, and 0.8, which are intended to segment the component into 5 layers of comparable size. (e) Temperature distributions along a vertical line in the middle of the component, with five different $\beta$ values.}
\label{fig:distribution_beta}
\end{figure}

The constant drain rate $\alpha_\texttt{T}$ is a key parameter in the heat equation. Alongside the thermal diffusivity $\kappa$, this parameter plays a significant role in determining the temperature distribution. To be able to generate layers of comparable size, the temperature distribution shall be smoothly varying between 0 and 1. Fig.~\ref{fig:distribution_beta} shows that a too small or large drain rate results in layers with significant variations in size. The domain is homogeneous in terms of thermal diffusivity $\kappa$, taking a value of $0.1$. The heat source, with a constant temperature of $1$, is placed on the bottom of the domain. If the drain rate is large, the temperature rapidly decreases to 0, as shown in Fig.~\ref{fig:distribution_beta} (d). After the transformation $\bm{t}=\bm{1}-\bm{\tau}$, and with uniformly sampled thresholds $\eta_j = \frac{j}{N}$, the layers are squeezed to the bottom. Reversely, a small drain rate also leads to uneven layers (Fig.~\ref{fig:distribution_beta} b). The segmentation in Fig.~\ref{fig:distribution_beta} (c) is more uniform, providing a good starting point for the optimization of the thermal diffusivity field towards the desired segmentation of the domain into layers of equal size. To ensure numerical convergence, the thermal diffusivity variable $\mu$ is constrained within the range of $0$ to $1$, matching the range of the density variable $\phi$. 

The temperature distribution is also affected by the size of the design domain. We introduce a size-independent parameter to help choose an appropriate drain rate. This parameter, denoted as $\beta$, relates the drain rate ($\alpha_\texttt{T}$) to the characteristic length of the design domain ($l_c$), 
\begin{equation}
    \alpha_\texttt{T} = \beta / l^{2}_{c}.
\label{eq:nondimensionalized}
\end{equation}
The unit of the drain rate is $\texttt{m}^{-2}$, while the unit of the length is $\texttt{m}$. Thus $\beta$ is non-dimensional. We choose the longest path of heat diffusion in the design domain as its characteristic length.

In Fig.~\ref{fig:distribution_beta} (e) we compare the temperature distribution under five different $\beta$ values. When $\beta = 1$, the resulting temperature smoothly decreases from 1 on the left, to close to 0 on the right. From numerical tests we observed good convergence in the optimization of $\mu$ and $\phi$ with $\beta \in [0.1, 1]$. Slight derivation from this range does not lead to dramatic convergence issues. It is robust for different domain sizes and volume fractions. In the results section, $\beta=0.1$ is consistently used for all examples.

\subsection{Formulation of space-time topology optimization}
To summarize, the space-time topology optimization using the heat equation is generically formulated as follows.
\begin{eqnarray}
\label{eq:formulation_thermal_density_sequence}
 & \underset{{\bm{\phi}, \; {\bm{\mu}}}}{\text{minimize}:} & f = f_{0}(\bm{\rho}) + \sum_{j=1}^{N}\alpha_jf_{j}(\bm{\rho}^{\{j\}}), \label{eq:obj}\\
 & \text{subject to}: &  \bm{K}(\bm{\rho})\bm{U} = \bm{F},  \label{eq:final_st}\\    
 &&  \bm{K}^{\{j\}}(\bm{\rho}^{\{j\}})\bm{U}^{\{j\}} = \bm{F}^{\{j\}}, \quad  j=1,...,N, \label{eq:inter_st} \\
     && \bm{K}_\texttt{T}(\bm{\phi}, \bm{\mu}) \bm{T} = \bm{b}, \label{eq:heatgovern}\\
     && g_0(\bm{\rho}) = \sum \rho_e v_e - V_0 \leq 0, \label{eq:volfrac}\\
 &&  g_j(\bm{\rho}^{\{j\}}) = \sum \rho_{e}^{\{j\}} v_e - \frac{j}{N}V_0 \le 0, \quad j=1,...,N, \label{eq:layer_vol_con}\\
  && 0 \leq \phi_e \leq 1, \forall e \in \mathcal{M},\\
 && 0 \leq \mu_e \leq 1, \forall e \in \mathcal{M}.
\end{eqnarray}
The objective function $f$ is composed of two parts. $f_{0}$ accounts for the performance of the final overall structure, while $f_{j}$ is related to each intermediate structure during the fabrication process. The different objectives are balanced by the weighting factors $\alpha_j$. The overall structure and intermediate structures all have their governing mechanical equation, i.e., Eq.~\ref{eq:final_st} and Eq.~\ref{eq:inter_st}, respectively. These mechanical equations are independent of each other, and thus can be solved in parallel. $\bm{K}$ and $\bm{U}$ are the global stiffness matrix and displacement vector corresponding to the overall structure, respectively, and $\bm{F}$ is the external force. The global stiffness matrix $\bm{K}^{\{j\}}$ of intermediate structures is determined by both the density variables $\phi$ and the diffusivity variables $\mu$. The load vector $\bm{F}^{\{j\}}$ may also be dependent on these variables, e.g., process-dependent thermomechanical loads. Eq.~\ref{eq:heatgovern} is the heat equation for obtaining a monotonic temperature field, as introduced in the previous subsection. The overall material usage is constrained in Eq.~\ref{eq:volfrac}, while an even amount of material for each intermediate structure is described by Eq.~\ref{eq:layer_vol_con} to ensure a constant fabrication speed.

Here, we present the general formulation of space-time topology optimization. For a more in-depth exploration of objective and constraint functions in specific application scenarios, such as structural compliance and thermal-induced distortion during fabrication, readers can refer to \cite{Wang2020SMO, Wang2023CMAME}. When considering the entire density field as fixed, the optimization model described above simplifies into the formulation outlined in Section~\ref{subsec:sqeuence_formulation}, focusing solely on sequence optimization. We solve the optimization problem using the method of moving asymptotes (MMA)~\cite{svanberg1987method}.

\subsection{Sensitivity analysis}
\label{subsec:sensitivity}

The general procedure for sensitivity calculation is consistent with those outlined in \cite{Wang2020SMO, Wang2023CMAME}. Here we focus on the primary distinction that arises from the fact that the time field is no longer a direct design variable; rather, it functions as an implicit variable dependent on both density $\phi$ and diffusivity $\mu$.

\paragraph{\textbf{Sensitivity w.r.t design variable $\bm{\mu}$}} Take the objective function $f_j$ as an example, the sensitivity is calculated using the adjoint method by adding an augmented Lagrangian term:
\begin{equation}
f_j = f_j +\bm{\lambda}^{\T}\left(\bm{K}_\texttt{T}\bm{T}-\bm{b}\right),
\end{equation}
where $\bm{\lambda}$ is the Lagrangian multiplier. Without loss of generality, if $f_j$ has no explicit component of the diffusivity design variable $\mu_e$, the derivative of $f_j$ with respect to $\mu_e$ is given by:
\begin{equation}
\begin{aligned}
\frac{\partial f_j}{\partial \mu_e}  &= \sum_{i \in \mathcal{M}} \frac{\partial f_j}{\partial t_i}\frac{\partial t_i}{\partial \mu_e} + \bm{\lambda}^{\T} \left(\frac{\partial \bm{K}_\texttt{T}}{\partial \mu_e}\bm{T} + \bm{K}_\texttt{T}\frac{\partial \bm{T}}{\partial \mu_e}\right).
\end{aligned}
\end{equation}
Rewrite the above equation in vector form while using the inverse relation between the time field and nodal temperature field $\bm{t}=\bm{1}-\bm{G}{\bm{T}}$:
\begin{equation}
\label{eq:adjoint_miu}
\begin{aligned}
\frac{\partial f_j}{\partial \mu_e} & = -\left(\frac{\partial f_j}{\partial \bm{t}}\right)^{\T} \bm{G} \frac{\partial \bm{T}}{\partial \mu_e} + \bm{\lambda}^{\T} \left(\frac{\partial \bm{K}_\texttt{T}}{\partial \mu_e}\bm{T} + \bm{K}_\texttt{T}\frac{\partial \bm{T}}{\partial \mu_e}\right) \\
& = \left(-\left(\frac{\partial f_j}{\partial \bm{t}}\right)^{\T} \bm{G} +\bm{\lambda}^{\T} \bm{K}_\texttt{T}\right)\frac{\partial \bm{T}}{\partial \mu_e} + \bm{\lambda}^{\T}\frac{\partial \bm{K}_\texttt{T}}{\partial \mu_e}\bm{T}
.
\end{aligned}
\end{equation}
To directly calculate the derivative of the temperature field with respect to the thermal diffusivity, the Lagrangian multiplier $\bm{\lambda}$ is chosen to eliminate the first term:
\begin{equation}
\label{eq:Lagrangian_miu}
-\left(\frac{\partial f_j}{\partial \bm{t}}\right)^{\T} \bm{G} +\bm{\lambda}^{\T} \bm{K}_\texttt{T} = 0.
\end{equation}
The derivative of $f_j$ regarding to $\mu_e$ can thus be rewritten as
\begin{equation}
\frac{\partial f_j}{\partial \mu_e} =\bm{\lambda}^{\T}\frac{\partial \bm{K}_\texttt{T}}{\partial \mu_e}\bm{T}.
\end{equation}

\paragraph{\textbf{Lagrangian multiplier $\bm{\lambda}$}} In order to solve for the Lagrangian multiplier in Eq.~\ref{eq:Lagrangian_miu}, the derivative of $f_j$ with respect to the time field $\bm{t}$ needs to be calculated. Since $f_j$ is the function of $\bm{\rho}^{\{j\}}$, from the chain rule:
\begin{equation}
\label{eq:chain_rule}
\frac{\partial f_j}{\partial t_e} =\frac{\partial f_j}{\partial \bar{t}_e}\frac{\partial \bar{t}_e}{\partial t_e}.
\end{equation}
where $\pdv*{\bar{t}_e}{t_e}$ is calculated from the Heaviside projection of Eq.~\ref{eq:projection}. In most of the cases, $f_j$ contains terms related to $\bm{U}^{\{j\}}$, thus another adjoint analysis is needed:
\begin{equation}
\begin{aligned}
\frac{\partial f_j}{\partial \bar{t}_e} & = \frac{\partial f_j}{\partial \bar{t}_e} +\bm{\gamma}^{\T}\pdv{\left(\bm{K}^{\{j\}}\bm{U}^{\{j\}}-\bm{F}^{\{j\}}\right)}{\bar{t}_e} \\
& =\frac{\partial f_j}{\partial \bm{U}^{\{j\}}}\frac{\partial \bm{U}^{\{j\}}}{\partial \bar{t}_e} +\bm{\gamma}^{\T}\left(\frac{\partial \bm{K}^{\{j\}}}{\partial \bar{t}_e}\bm{U}^{\{j\}} + \bm{K}^{\{j\}}\frac{\partial \bm{U}^{\{j\}}}{\partial \bar{t}_e}\right) \\
& = \left( \frac{\partial f_j}{\partial \bm{U}^{\{j\}}} + \bm{\gamma}^{\T}  \bm{K}^{\{j\}}\right)\frac{\partial \bm{U}^{\{j\}}}{\partial \bar{t}_e} + \bm{\gamma}^{\T}\frac{\partial \bm{K}^{\{j\}}}{\partial \bar{t}_e}\bm{U}^{\{j\}}.
\end{aligned}
\end{equation}
Here $\bm{\gamma}$ is the Lagrangian multiplier related to the mechanical equilibrium equation of intermediate structure, which is solved from:
\begin{equation}
\label{eq:lagrangian_mechanical}
\frac{\partial f_j}{\partial \bm{U}^{\{j\}}} + \bm{\gamma}^{\T}  \bm{K}^{\{j\}} = 0.
\end{equation}
The equation above simply assume that there is no explicit term about $\bar{t}_e$ in $f_j$, and the external loads $F^{\{j\}}$ is independent on the design $\bm{\rho}^{\{j\}}$. These are not necessary true, for example, design-dependent load like self-weight can also be considered by making slight modifications on the sensitivity analysis.

\paragraph{\textbf{Sensitivity w.r.t design variable $\bm{\phi}$}} The derivative of $f_j$ with respect to the density design variable $\phi_e$ can be obtained using the chain rule, which aligns with the filtering and projection procedure:
\begin{equation}
\label{eq:chain_rule}
\frac{\partial f_j}{\partial \phi_e} =\sum_{i\in \mathcal{N}_e} \dv{f_j}{\rho_i}\frac{\partial \rho_i}{\partial \tilde{\phi}_i}\frac{\partial \tilde{\phi}_i}{\partial \phi_e}.
\end{equation}
The derivative of $f_j$ regarding the physical density $\rho_e$ is also calculated using the adjoint method. However, it should be noted that since the objective function has multiple components that contains $\rho_e$ implicitly. As indicated in Eq.~\ref{eq:heatconductivity}, the time field obtained from the heat equation is also dependent on the density field. We thus divide the sensitivity of $f_j$ into two parts. The first part is independent of the time field $\bm{t}$, while the second part is related to the time field $\bm{t}$. 

\begin{equation}
\begin{aligned}
\dv{f_j}{\rho_e}  &= \frac{\partial f_j}{\partial \rho_e} + \sum_{i \in \mathcal{M}} \frac{\partial f_j}{\partial t_i}\frac{\partial t_i}{\partial \rho_e}.
\end{aligned}
\end{equation}
Similar to Eq.~\ref{eq:adjoint_miu}, the above equation is written into an augmented form:
\begin{equation}
\label{eq:adjoint_rho}
\begin{aligned}
\dv{f_j}{\rho_e} & = \frac{\partial f_j}{\partial \rho_e} -\left(\frac{\partial f_j}{\partial \bm{t}}\right)^{\T} \bm{G} \frac{\partial \bm{T}}{\partial \rho_e} + \bm{\lambda}^{\T} \left(\frac{\partial \bm{K}_\texttt{T}}{\partial \rho_e}\bm{T} + \bm{K}_\texttt{T}\frac{\partial \bm{T}}{\partial \rho_e}\right) \\
& = \frac{\partial f_j}{\partial \rho_e} + \left(-\left(\frac{\partial f_j}{\partial \bm{t}}\right)^{\T} \bm{G} +\bm{\lambda}^{\T} \bm{K}_\texttt{T}\right)\frac{\partial \bm{T}}{\partial \rho_e} + \bm{\lambda}^{\T}\frac{\partial \bm{K}_\texttt{T}}{\partial \rho_e}\bm{T}
.
\end{aligned}
\end{equation}
By using the same Lagrangian multiplier as solved from Eq.~\ref{eq:Lagrangian_miu}, the above equation is simplified into
\begin{equation}
\begin{aligned}
\dv{f_j}{\rho_e}  &= \frac{\partial f_j}{\partial \rho_e} + \bm{\lambda}^{\T}\frac{\partial \bm{K}_\texttt{T}}{\partial \rho_e}\bm{T}.
\end{aligned}
\end{equation}
The first part is directly related to $\phi_e$, and can also be calculated by using adjoint analysis:
\begin{equation}
\begin{aligned}
\frac{\partial f_j}{\partial \rho_e}  &= \frac{\partial f_j}{\partial \rho_e} + \bm{\gamma}^{\T}\pdv{\left(\bm{K}^{\{j\}}\bm{U}^{\{j\}}-\bm{F}^{\{j\}}\right)}{\rho_e} \\
&= \left( \frac{\partial f_j}{\partial \bm{U}^{\{j\}}} + \bm{\gamma}^{\T}  \bm{K}^{\{j\}}\right)\frac{\partial \bm{U}^{\{j\}}}{\partial \rho_e} + \bm{\gamma}^{\T}\frac{\partial \bm{K}^{\{j\}}}{\partial \rho_e}\bm{U}^{\{j\}}\\
&= \bm{\gamma}^{\T}\frac{\partial \bm{K}^{\{j\}}}{\partial \rho_e}\bm{U}^{\{j\}}.
\end{aligned}
\end{equation}
Here $\bm{\gamma}$ is same the Lagrangian multiplier as in  Eq.~\ref{eq:lagrangian_mechanical}.

The sensitivity analysis of these volume constraint functions $g_0$ and $g_j$ is a simpler version of the analysis above, as they only have explicit parts related to their corresponding $\bm{\rho}^{\{j\}}$.


\section{Results and discussion}
\label{sec:results}

In this section, we first demonstrate that the effectiveness of the proposed method is independent of the initialization of the optimization variables, in contrast to prior formulations with explicit constraints (Section~\ref{subsec:validate}). We then evaluate our method for different application scenarios, including the minimization of gravity-induced distortion of intermediate structures during construction (Section~\ref{subsec:gravity}), and  minimization of thermal-induced distortion of fixed structural layout (Section~\ref{subsec:thermal_component}) and by concurrent optimization of the structural layout and fabrication sequence (Section~\ref{subsec:thermal_density}). Lastly, we test the proposed method for space-time topology optimization of 3D structures (Section~\ref{subsec:3D}). In all experiments, the Young's modulus was fixed at 1.0, and the Poisson's ratio was fixed at 0.3. Unless stated otherwise, in all 2D tests, we set 20 layers and a volume fraction of 0.6. The method was implemented in Matlab.

\subsection{Validation of avoiding local minima}
\label{subsec:validate}

\begin{figure*}[htb!]
\centering
\includegraphics[width=.7\linewidth]{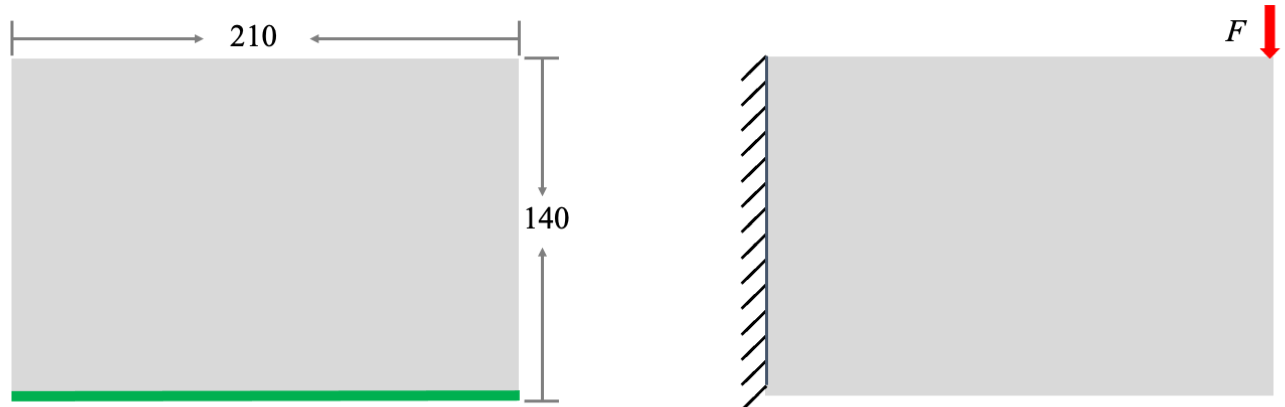}
\caption{Illustration of the problem setup. Left: the build plate (green) is placed at the bottom of the design domain, meaning the structure is to be constructed from the bottom upwards. Right: the boundary condition of the structure in its intended use scenario.}
\label{fig:design-domain}
\end{figure*}

The first example is intended to validate that our new formulation avoids the appearance of local minima in the time field. The test is performed on a cantilever beam, as illustrated in Fig.~\ref{fig:design-domain}. The design domain is discretized into a regular quadrilateral mesh with the resolution of $210\times140$. The objective is to reduce the compliance of the entire component, under a load on the top-right corner. Process-dependent loads on intermediate structures are not considered in this test, by setting $\alpha_j$ in Eq.~\ref{eq:obj} to zero. To test the robustness of the new formulation, we run space-time topology optimization with different initializations of the optimization variables. We study in particular different initializations of the thermal diffusivity field $\bm{\mu}$. For the density-related field $\bm{\phi}$, we initialize it with a uniform density of the target volume fraction. 

\begin{figure*}[htb!]
\centering
\includegraphics[width=0.98\linewidth]{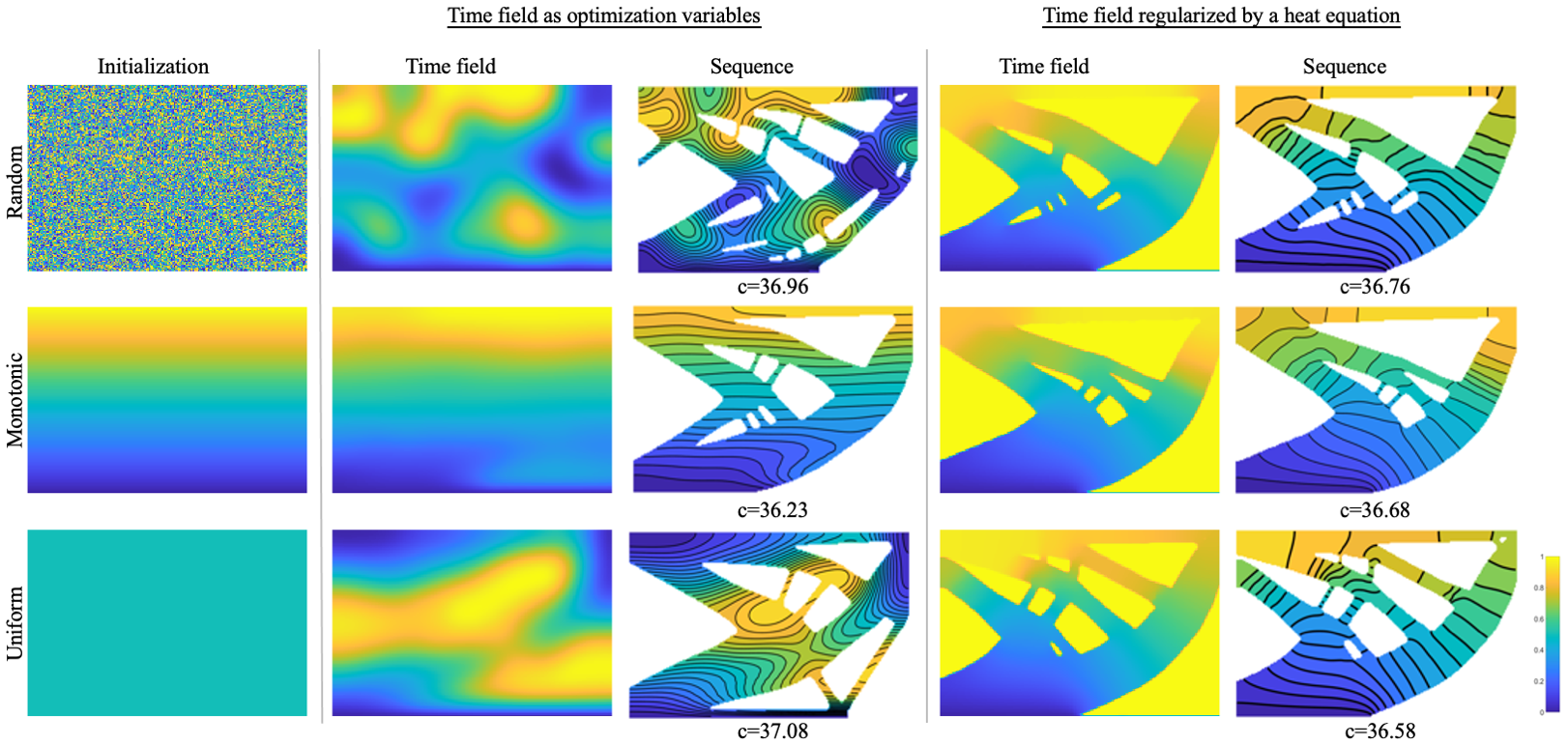}
\caption{Comparison of space-time topology optimization with the time field as optimization variables~\cite{Wang2020SMO} and using the proposed regularization. The three different initializations are, from top to bottom, random, monotonic, and uniform. When using the time field as optimization variables, the time field is initialized, while in our proposed method, the initialization is applied to the thermal diffusivity field $\bm{\mu}$. Using the proposed regularization, the optimized time field is monotonic, regardless of the initial values.}
\label{fig:initialization}
\end{figure*}

Figure~\ref{fig:initialization} compares the results obtained using the original formulation~\cite{Wang2020SMO} and the proposed new formulation, under different initializations. All results were generated with 500 optimization iterations. The three initializations are shown in the first column,  from top to bottom, random, monotonic, and uniform. The second and third columns visualize the optimized time field and corresponding fabrication sequence, respectively, obtained using the original formulation where the optimization variables are the time field~\cite{Wang2020SMO}. Under random as well as uniform initialization of the time field, the optimized time fields exhibit multiple local minima and maximum. Only the monotonic initialization results in a clear fabrication sequence. The last two columns visualize the optimized results with the proposed regularization. Here, regardless of the initialization of the thermal diffusivity field $\bm{\mu}$, the optimized time fields are monotonic. The time values increase, within the shape of the optimized structural designs, from the bottom which is prescribed as the build plate, to the up boundary of the domain. The areas outside of the optimized structural designs take a time value of $1$, i.e., the upper bound of time values. The monotonic increase of the time field is also clearly visible from the optimized sequence shown on the right-hand side. In all tests, the compliance values are comparable. This is because the time field in this example is only used to segment the evolving structural layout into layers of equal areas, and no sequence-dependent loads are applied. 

\begin{figure*}[htb!]
\centering
\includegraphics[width=0.98\linewidth]{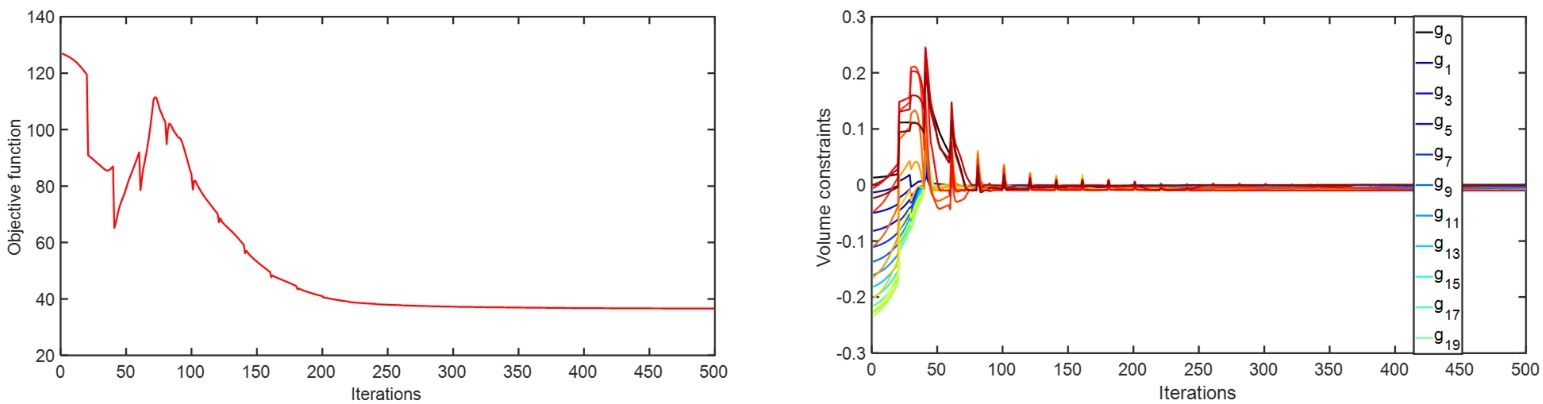}
\caption{The convergence curves of space-time optimization using the proposed regularization. Left: the objective function. Right: the per-layer volume constraints.}
\label{fig:converge}
\end{figure*}

Fig.~\ref{fig:converge} plots the convergence curves of the objective function (left) and volume constraint per layer (right), corresponding to the optimization with a uniform thermal diffusivity field. The optimization stabilizes after about 200 iterations. A continuation is applied to the projection parameter $\beta_t$ in Eq.~\ref{eq:projection}. From these curves, we can see, that the objective function converges well and all volume constraints are satisfied.

\subsection{Space-time topology optimization under gravity loads}
\label{subsec:gravity}

Our next example is the concurrent optimization of the structural design and fabrication sequence under gravity loads, a problem that has been studied in Ref.~\cite{Wang2020SMO}. During fabrication, the self-weight of the intermediate structure may result in large deflection. The weight of the intermediate structure increases as the fabrication progresses, and it is dependent on the fabrication sequence. In this example, the objective function consists of two parts. The first part is the compliance of the final structure, under the boundary condition as illustrated in Fig.~\ref{fig:design-domain}(right). The second part is the sum of the compliance values of intermediate structures under self-weight. These two parts are integrated with a weighting factor $\alpha$ in Eq.~\ref{eq:obj}. We note that the boundary condition for simulating the effects of gravity on the intermediate structure is different than that for simulating the load-bearing capacity of the complete structure.

\begin{figure*}[htb!]
\centering
\includegraphics[width=0.98\linewidth]{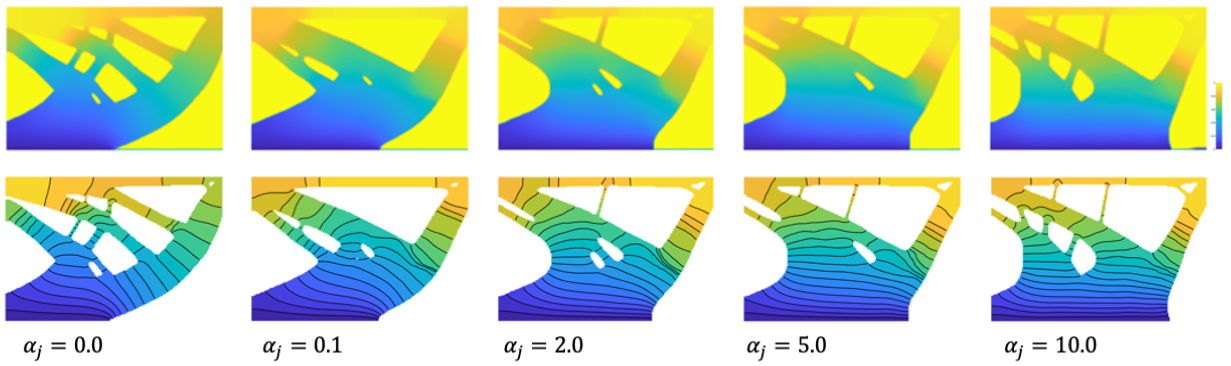}
\caption{Space-time topology optimization with gravity loads on intermediate structures. From left to right, the weighting factor for intermediate structures increases from $0.0$ to $10.0$. Top: the optimized time fields. Bottom: optimized structural layouts with the corresponding fabrication sequence.}
\label{fig:gravity}
\end{figure*}

Fig.~\ref{fig:gravity} shows the optimization results with five different weighting factors. From left to right, $0.0$, $0.1$, $2.0$, $5.0$ and $10.0$. The first row shows the optimized time fields. The second row shows the optimized structure as well as the fabrication sequence according to the time field. In all five cases, the fabrication sequence starts from the bottom of the domain where the build plate is located, and gradually moves upwards. Comparing the five results from left to right, as the weighting factor for intermediate structures increases, the optimized structure exhibits a progressively larger contact region with the build plate and concentrates more material in its lower part. 

\subsection{Fabrication sequence optimization for minimizing thermal-induced distortion}
\label{subsec:thermal_component}

In the third set of experiments, we validate our method in fabrication sequence optimization for minimizing the distortion induced by process-dependent thermomechanical loads. Distortion is a major issue in metal additive manufacturing. It is caused by the high energy input and large temperature gradients acting on the materials during fabrication. In a prior work~\cite{Wang2023CMAME}, the time field was directly optimized to obtain the fabrication sequence. As shown in Fig.~\ref{fig:distance}, the feasibility of the optimized sequence depends on the initial values of the time field.

\begin{figure*}[htb!]
\centering
\includegraphics[width=0.9\linewidth]{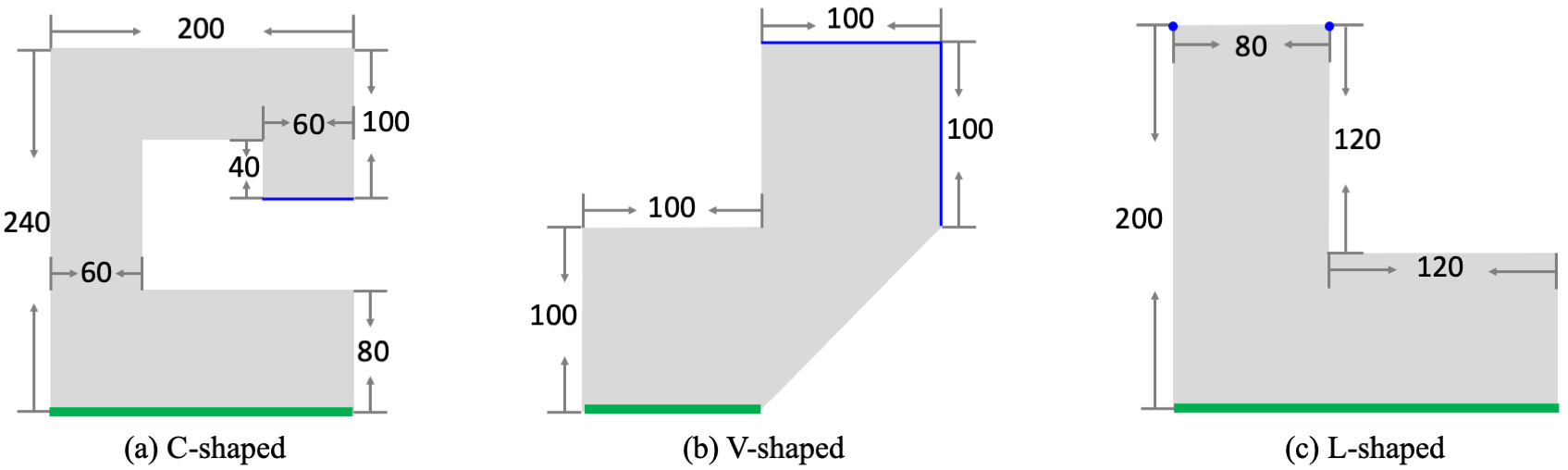}
\caption{Illustration of the three components used in fabrication sequence optimization. The build plate (green) is located at the bottom of each domain. The objective is to minimize the thermal-induced distortion, measured on the geometric features marked in blue.}
\label{fig:components}
\end{figure*}

We test the proposed regularization on sequence optimization for three components, shown in Fig.~\ref{fig:components}. The green line at the bottom of each component indicates the build plate. The objective is to minimize thermal-induced distortion. In (a) the distortion is represented by the flatness of the blue edge, measured by the relative vertical displacements of nodes on the edge. In (b) the distortion is represented by the flatness of a horizontal and a vertical edge, both indicated by blue lines. In (c) the distortion is represented by the relative vertical displacements of two nodes on the top edge, indicated by blue dots. The distortion due to process-dependent thermomechanical loads is simulated using the inherent strain method. We refer readers to \cite{Wang2023CMAME} for a detailed explanation of the simulation of the layered manufacturing process.

\begin{figure*}[htb!]
\centering
\includegraphics[width=.9\linewidth]{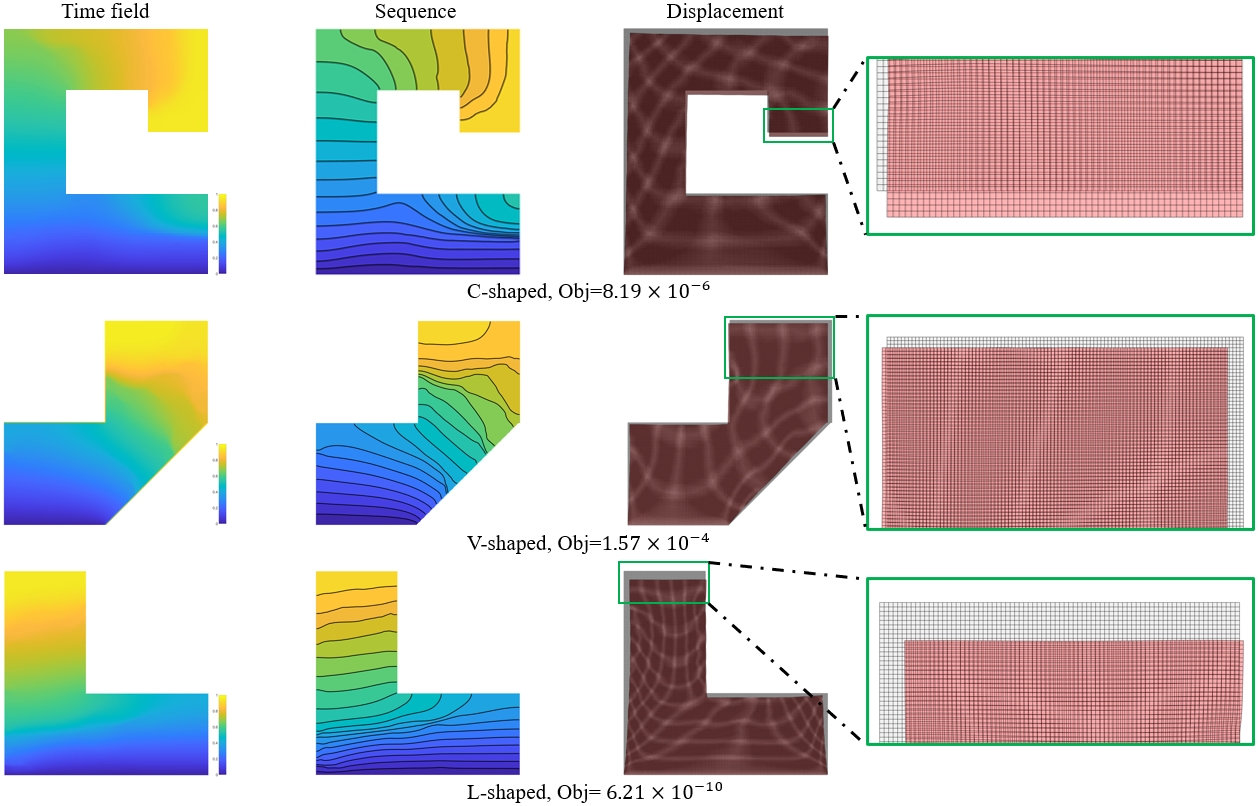}
\caption{Results of fabrication sequence optimization for distortion minimization. From left to right: the optimized time field, fabrication sequence, and simulated displacement field.}
\label{fig:therm_given_component}
\end{figure*}

The optimized results are visualized in Fig.~\ref{fig:therm_given_component}. From left to right, the optimized time field, corresponding fabrication sequence, and resulting distortion. The optimized time field of each component shows a continuous increase from the bottom of the domain to an edge that is furthest. This continuous change can also be observed in the optimized fabrication sequence. In all three examples, the measured distortion values are very small. From the enlarged figures on the right, it can be seen that the intended flatness of the edges is effectively achieved. 

\subsection{Concurrent structure and sequence optimization with constraints on thermal-induced distortion}
\label{subsec:thermal_density}

\begin{figure*}[htb!]
\centering
\includegraphics[width=0.9\linewidth]{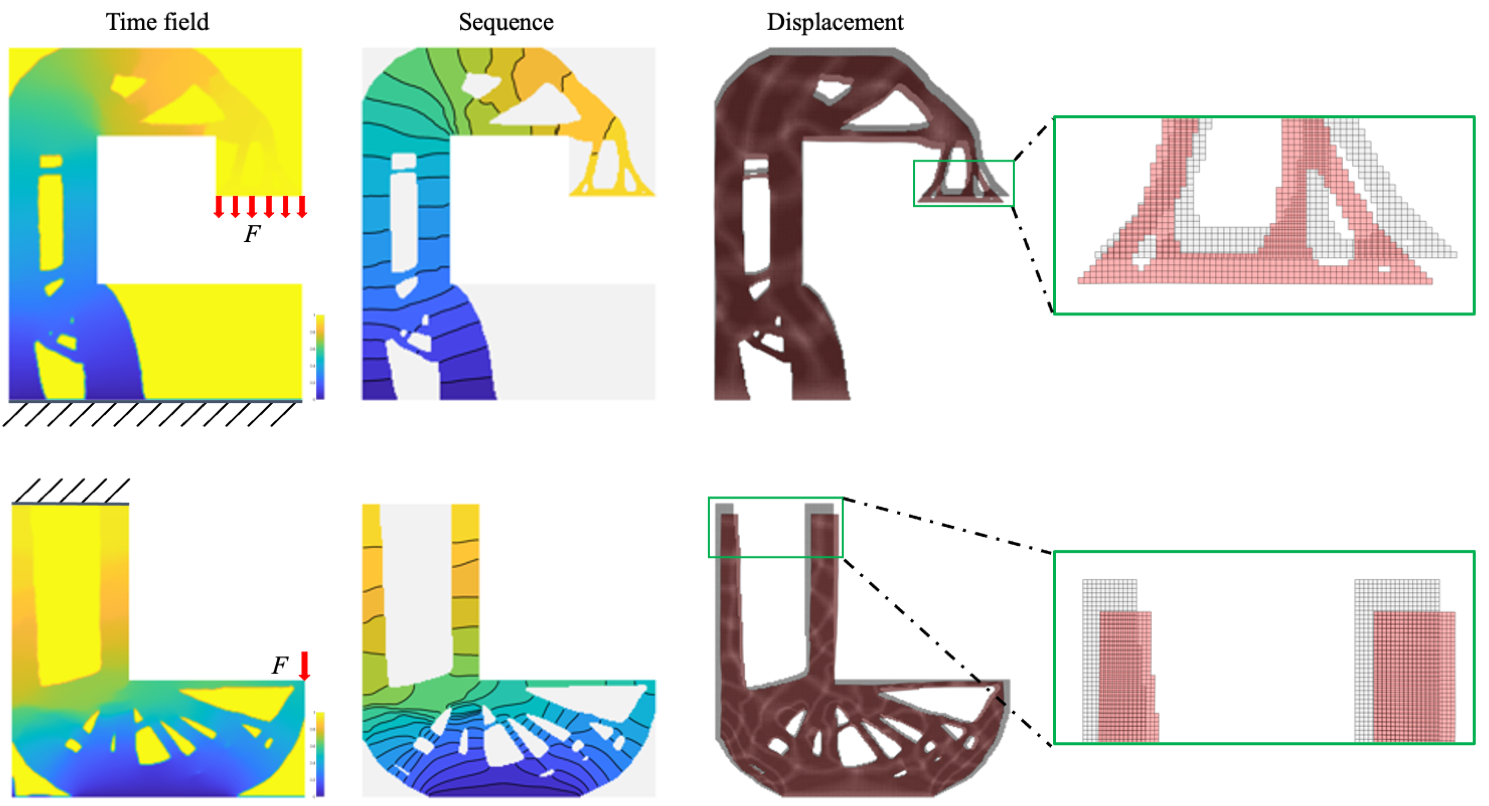}
\caption{Results of space-time topology optimization under constraints on thermal-induced distortion. From left to right: the optimized time field, fabrication sequence, and simulated displacement field.}
\label{fig:density_sequence}
\end{figure*}

To further demonstrate the effectiveness of the proposed method, we concurrently optimize the structural layout and fabrication sequence. The test was performed on the design domains shown in Fig.~\ref{fig:therm_given_component}(a) and (c). The objective is to minimize the compliance of the final structure, under the boundary condition illustrated in Fig.~\ref{fig:density_sequence}(left). The optimization is subject to a bound on the thermal-induced distortion. The measurement of distortion follows the previous subsection. The target volume fractions of the two components are 0.5 and 0.6, respectively.

Fig.~\ref{fig:density_sequence} visualizes the optimized time fields, sequences, and displacements. While the optimized structural layouts are very complex, the time fields are monotonic within the solid regions, demonstrating the effectiveness of the proposed method.

\subsection{Space-time topology optimization of 3D structures}
\label{subsec:3D}

\begin{figure*}[htb!]
\centering
\includegraphics[width=0.98\linewidth]{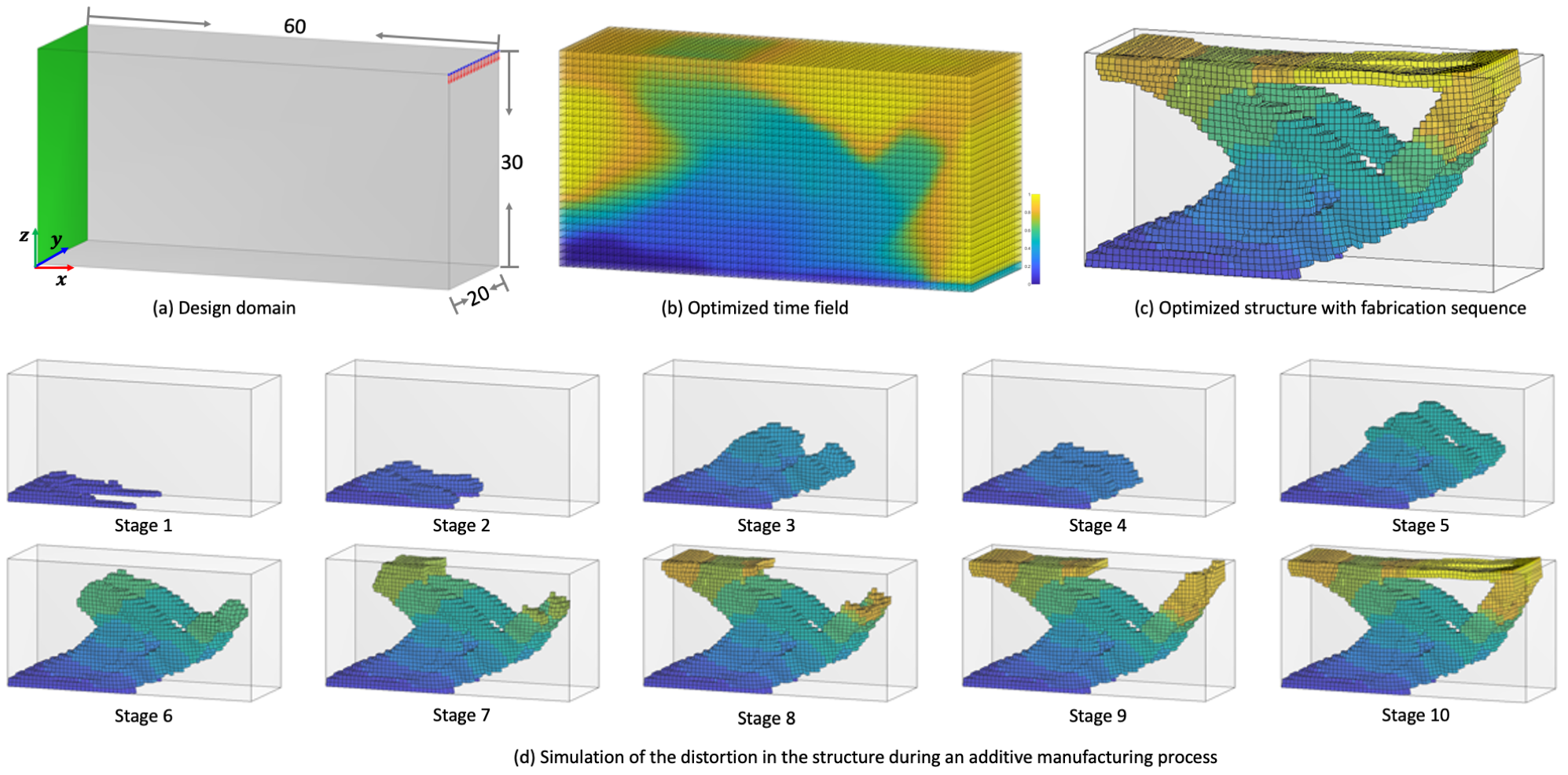}
\caption{Space-time topology optimization of a 3D structure. (a): the design domain and boundary condition. (b): the optimized time field. (c): simulation of the distortion in the optimized structural layout with the corresponding fabrication sequence. (d): Simulation of intermediate structures during the fabrication process.} 
\label{fig:3D_optimized}
\end{figure*}

The proposed regularization method is also tested on space-time topology optimization of 3D structures. The design domain and boundary condition are illustrated in Fig.~\ref{fig:3D_optimized} (a). The left side of the domain is fixed, while a distributed load is applied on the top-right edge. The objective is to minimize compliance under this boundary condition. This is subject to a constraint on thermal-induced distortion, which is represented by the flatness of the edge where the external load is applied, i.e., the relative displacement of nodes on this edge along the $x$- and $z$-axis due to process-dependent thermomechanical loads. The build plate is situated at the bottom of the domain, meaning the construction of the structure will proceed from the bottom up. The design domain is discretized into hexahedral finite elements with the resolution of $60\times20\times30$. The volume fraction is 0.2, and the number of layers is 10. 

The optimized time field is shown in Fig.~\ref{fig:3D_optimized} (b), while the optimized structure under thermal-induced distortion is visualized in (c). The flatness of the target edge is clearly visible. Also shown in (c) are the segmented layers, indicated by the varying colors. A simulation of the fabrication sequence is shown in (d).

\section{Conclusions}
\label{sec:conclusions}
In this paper, we present a novel method for regularizing the time field in space-time topology optimization. Instead of directly treating the time field as optimization variables, we employ a static heat equation with a drain term to derive a monotonic pseudo-temperature field. We then consider spatially varying thermal diffusivity in the heat equation as optimization variables. This approach effectively prevents local minima from appearing in the optimized time field, regardless of the initial values used in optimization. This approach is easily extendable to deal with 3D components. We demonstrate its effectiveness in additive manufacturing scenarios, which involve process-dependent loads such as gravity and thermomechanical loads on intermediate structures.

A monotonically progressing time field is fundamental for ensuring the viability of the interpreted fabrication sequence. This work thus lays a solid foundation for further advancements in space-time optimization for multi-axis additive manufacturing. As part of our future work, we are exploring methods to regularize the thickness of curved layers. While multi-axis additive manufacturing facilitates curved layers, achieving uniform thickness, or near uniformity, remains a desirable goal.

\vspace{10mm}
\small{\textbf{Acknowledgements} The authors gratefully acknowledge the support from the LEaDing Fellows Programme at the Delft University of Technology, which has received funding from the European Union's Horizon 2020 research and innovation programme under the Marie Skłodowska-Curie grant agreement No. 707404. Weiming Wang wishes to thank the National Natural Science Foundation of China (No. 62172073). Kai Wu is supported by the China Scholarship Council (CSC). This work is also partly supported by the Dutch Research Council (NWO) with project number 20382, ``Space-time optimization for additive manufacturing''.
}

\bibliography{mybib}

\end{document}